\documentclass{amsart}
\usepackage{amsfonts}
\usepackage{color}
\usepackage{hyperref}
\usepackage{cite}

\newtheorem{theorem}{Theorem}
\theoremstyle{plain}

\newtheorem{corollary}{Corollary}

\newtheorem{lemma}{Lemma}

\newtheorem{proposition}{Proposition}

\numberwithin{equation}{section}

\begin{document}
\title[Concircular curvature Tensor on Warped Products]{Concircular
curvature on warped product manifolds and applications}
\author{Uday Chand De}
\address[U. C. De]{Department of Pure Mathematics, University of Calcutta,
35 Ballygaunge Circular Road, Kolkata 700019, West Bengala, India}
\email{uc$\_$de@yahoo.com}
\author{Sameh Shenawy}
\address[S. Shenawy]{Basic Science Department, Modern Academy for
Engineering and Technology, Maadi, Egypt,}
\email{drssshenawy@eng.modern-academy.edu.eg, drshenawy@mail.com}
\author{B\"ulent \"Unal}
\address[B. \"{U}nal]{Department of Mathematics, Bilkent University,
Bilkent, 06800 Ankara, Turkey}
\email{bulentunal@mail.com}

\begin{abstract}
This study aims mainly at investigating the effects of concircular flatness
and concircular symmetry of a warped product manifold on its fibre and base
manifolds. Concircularly flat and concircularly symmetric warped product
manifolds are investigated. The divergence free concircular curvature tensor
on warped product manifolds is considered. Finally, we apply some of these
results to generalized Robertson-Walker and standard static space-times.
\end{abstract}

\maketitle

\section{Introduction}

A transformation which preserves geodesic circles is called a concircular
transformation\cite{Yano:1940}. The geometry which deals with concircular
transformation is called concircular geometry. The concircular curvature
tensor $\mathcal{C}$ remains invariant under concircular transformation of a
(pseudo-)Riemannian manifold $M$. $M$ is called concircularly flat if its
concircular curvature tensor $\mathcal{C}$\ vanishes at every point. A
concircularly flat manifold $M$ is a manifold of constant curvature. Thus
the tensor $\mathcal{C}$ measures the deviation of $M$ from constant
curvature (For further details, see \cite{Ahsan:2017,Yano:1940}).

In a series of studies, Pokhariyal and Mishra studied the recurrent
properties and relativistic significance of concircular curvature tensor,
amongst many others, in Riemannian manifolds\cite{Mishra:1971,
Pokhariyal:1970,Pokhariyal:1974,Pokhariyal:1982}. Concircularly
semi-symmetric $K-$contact manifolds are considered in \cite{Majhi:2013} and
concircularly recurrent Finsler manifolds are studied in \cite{Youssef:2013}%
. In \cite{Majhi:2015}, the authors considered $N(k)-$contact metric
manifolds satisfying $\mathcal{C}\cdot \mathcal{P}=0$, where $\mathcal{P}$
denotes the projective curvature tensor. Similarly, a study of $\left( k,\mu
,\nu \right) -$contact metric $3-$manifolds satisfying one of the conditions 
$\nabla \mathcal{C}=0$, $\mathcal{C}\left( \zeta ,X\right) \cdot \mathcal{C}%
=0$, $\mathcal{R}\left( \zeta ,X\right) \cdot \mathcal{C}=0,$ where $\zeta $
is the Reeb field, is considered in \cite{Gouli:2014}. Perfect fluid
space-times with either vanishing or divergence-free concircular curvature
tensor are considered in \cite{Ahsan:2009}. The authors of \cite%
{Zlatanovic:2014} considered equitorsion concircular mapping between
generalized Riemannian manifolds (in the sense of Eisenhart's definition)
and obtained some invariant curvature tensors. These tensors are
generalizations of concircular curvature tensor on Riemannian manifolds. In
a recent paper \cite{Chen2015}, Chen provided some classification of Ricci
solitons with respect to a concircular potential field. In \cite{Kim1982},
the concept of special concircular vector fields is introduced and it is
proved that an $n$-dimensional Riemannian manifold that admits $n$ linearly
independent special concircular vector fields has constant sectional
curvature. Similarly, in \cite{Catino2105}, the authors characterize the
local structure of a Riemannian manifold whose Codazzi tensor has exactly
two distinct eigenvalues. In \cite{Olszak:2012}, it is proven that each
concircularly recurrent manifold is necessarily a recurrent manifold.

Motivated by these studies and many others, the main purpose of this article
is to study concircular curvature tensor on warped product manifolds and to
apply some of the results to two different $n-$dimensional space-times,
namely, generalized Robertson-Walker space-times and standard static
space-times. Concircularly flat and concircularly symmetric warped product
manifolds are also considered. Finally, divergence free concircular
curvature tensor on warped product manifolds is investigated.

This article is organized as follows. The next section presents the main
properties of the concircular curvature tensor. In Section 3, the
semi-symmetries of the concircular curvature tensor are investigated.
Section 4 is devoted to the study of concircularly flat warped product
manifolds whereas Section 5 is devoted to the study of concircularly
symmetric warped product manifolds. Finally, divergence free concircular
curvature tensor on warped product space-time models is considered in
Section 6.

\section{Concircular Curvature Tensor}

Let $\left( M,g\right) $ be a pseudo-Riemannian $n-$dimensional manifold $%
n\geq 3$. Throughout this section, $\nabla ,$ $\mathcal{R}$, {\textrm{Ric}}
and $\tau $ denote the Levi-Civita connection, curvature tensor, Ricci
curvature and scalar curvature of the metric tensor $g,$ respectively.

The concircular curvature tensor $\mathcal{C}$ on a pseudo-Riemannian
manifold $\left( M,g,\nabla \right) $ is defined as follows \cite%
{Ahsan:2009,Pokhariyal:1970, Singh:1999,Yano:1984}. Let $X,Y,Z,V\in 
\mathfrak{X}\left( M\right) $, then%
\begin{eqnarray}
\mathcal{C}\left( X,Y\right) Z &=&\mathcal{R}\left( X,Y\right) Z  \notag \\
&&-\frac{\tau }{n\left( n-1\right) }\left[ g\left( X,Z\right) Y-g\left(
Y,Z\right) X\right] ,
\end{eqnarray}%
where $\mathcal{R}\left( X,Y\right) Z=\nabla _{Y}\nabla _{X}Z-\nabla
_{X}\nabla _{Y}Z+\nabla _{\left[ X,Y\right] }Z$ is the Riemann curvature
tensor. It is clear that $\mathcal{C}\left( X,Y\right) Z$ is skew-symmetric
in the first two indices. Furthermore, 
\begin{eqnarray}
\mathcal{C}\left( X,Y,Z,V\right) &=&\mathcal{R}\left( X,Y,Z,V\right)  \notag
\\
&&-\frac{\tau }{n\left( n-1\right) }\left[ g\left( X,Z\right) g\left(
Y,V\right) -g\left( Y,Z\right) g\left( X,V\right) \right] .  \label{E0}
\end{eqnarray}

The definition of the concircular curvature tensor in local coordinates is
as follows%
\begin{equation}
\mathcal{C}_{ijkl}=\mathcal{R}_{ijkl}-\frac{\tau }{n\left( n-1\right) }\left[
g_{ik}g_{jl}-g_{jk}g_{il}\right] ,
\end{equation}%
where $\tau =g^{ij}\mathcal{R}_{ij}$ is the scalar curvature. This formula
suggests a generalization of this tensor of the form%
\begin{equation}
\mathcal{K}_{ijkl}=a_{0}\mathcal{R}%
_{ijkl}+a_{1}g_{ij}g_{kl}+a_{2}g_{ik}g_{jl}+a_{3}g_{jk}g_{il},  \label{E1}
\end{equation}%
where $a_{i}$ are constants and $a_{0}\neq 0$. Assume that a
pseudo-Riemannian manifold $\left( M,g\right) $ is a $\mathcal{K-}$curvature
flat manifold, then%
\begin{equation}
a_{0}\mathcal{R}%
_{ijkl}+a_{1}g_{ij}g_{kl}+a_{2}g_{ik}g_{jl}+a_{3}g_{jk}g_{il}=0.  \label{E2}
\end{equation}%
Multiplying both sides by $g^{il},$ we get%
\begin{equation*}
-a_{0}\mathcal{R}_{jk}+a_{1}g_{kj}+a_{2}g_{kj}+na_{3}g_{jk}=0.
\end{equation*}%
Again, by multiplying both sides of Equation (\ref{E2}) by $g^{ik}$, we get%
\begin{equation*}
a_{0}\mathcal{R}_{jl}+a_{1}g_{jl}+na_{2}g_{jl}+a_{3}g_{jl}=0
\end{equation*}%
Thus $\left( M,g\right) $ is Einstein with%
\begin{eqnarray*}
\mathcal{R}_{jk} &=&\frac{a_{1}+a_{2}+na_{3}}{a_{0}}g_{jk} \\
\mathcal{R}_{jl} &=&-\left( \frac{a_{1}+na_{2}+a_{3}}{a_{0}}\right) g_{jl}.
\end{eqnarray*}%
Consequently, a second contraction implies%
\begin{eqnarray*}
a_{1}+a_{2}+na_{3} &=&\frac{\tau }{n}a_{0} \\
a_{1}+na_{2}+a_{3} &=&-\frac{\tau }{n}a_{0}.
\end{eqnarray*}%
However, Equation (\ref{E2}) yields%
\begin{equation*}
na_{1}+a_{2}+a_{3}=0.
\end{equation*}%
These equations imply that%
\begin{equation*}
a_{1}=0,\ \ \ a_{2}=-a_{3}=\frac{-a_{0}\tau }{n\left( n-1\right) }
\end{equation*}%
Again, Equation (\ref{E2}) becomes%
\begin{equation*}
a_{0}\left[ \mathcal{R}_{ijkl}-\frac{\tau }{n\left( n-1\right) }\left(
g_{ik}g_{jl}-g_{jk}g_{il}\right) \right] =0.
\end{equation*}

Thus, $M$ is of constant sectional curvature. Therefore, the only $\mathcal{%
K-}$curvature flat tensor is the concircular curvature tensor and we have:

\begin{theorem}
Let $M$ be a $\mathcal{K-}$curvature flat manifold. Then $M$ is of constant
curvature and%
\begin{equation}
a_{1}=0,\ \ \ a_{2}=-a_{3}=\frac{-a_{0}\tau }{n\left( n-1\right) }
\label{E3}
\end{equation}%
i.e. $\mathcal{K}$ is a constant multiple of $\mathcal{C}$.
\end{theorem}

This result is proved in \cite{Yano:1940} when $\mathcal{K=C}$. Moreover, it
is found in \cite{Youssef:2013} for the Finslerian case. Assume that $M$ be
a $4-$dimensional space-time obeying Einstein's field equation with
cosmological constant, i.e.%
\begin{equation}
\mathcal{R}_{ij}-\frac{\tau }{2}g_{ij}+\Lambda g_{ij}=kT_{ij}  \label{E4}
\end{equation}%
where $T$ is the energy-momentum tensor. Let us define%
\begin{eqnarray}
\mathcal{K}_{jl} &=&g^{ik}\mathcal{K}_{ijkl}  \notag \\
&=&a_{0}\mathcal{R}_{jl}+\left( a_{1}+na_{2}+a_{3}\right) g_{jl}.  \label{E5}
\end{eqnarray}%
Now Equation (\ref{E4}) becomes%
\begin{equation}
\mathcal{K}_{ij}-\left( a_{1}+na_{2}+a_{3}\right) g_{ij}-\frac{a_{0}\tau }{2}%
g_{ij}+a_{0}\Lambda g_{ij}=a_{0}kT_{ij}.  \label{E6}
\end{equation}%
Thus%
\begin{equation*}
\nabla _{i}\mathcal{K}_{j}^{i}-\frac{a_{0}}{2}\nabla _{j}\tau =a_{0}k\nabla
_{i}T_{j}^{i}.
\end{equation*}%
We can now state:

\begin{theorem}
In a relativistic space-time obeying Einstein's field equations, the
energy-momentum tensor is divergence free if and only if%
\begin{equation*}
\nabla _{i}\mathcal{K}_{j}^{i}=\frac{a_{0}}{2}\nabla _{j}\tau .
\end{equation*}
\end{theorem}

\section{Semi-symmetries of $\mathcal{C}$}

It is noted that Equation (\ref{E0}) has the form%
\begin{eqnarray}
\mathcal{C} &=&\mathcal{R}-\frac{\tau }{2n\left( n-1\right) }\left( g\wedge
g\right)   \notag \\
&=&\mathcal{R}-\frac{\tau }{n\left( n-1\right) }\mathcal{G},
\end{eqnarray}%
where $\wedge $ is Kulkarni$-$Nomizu product of two symmetric $2-$tensors
(see \cite[p. 47]{Besse:2008}) and $\mathcal{G}=\frac{1}{2}\left( g\wedge
g\right) $. This equation leads us to%
\begin{eqnarray*}
\mathcal{R\cdot C} &=&\mathcal{R\cdot R}-\frac{\tau }{n\left( n-1\right) }%
\mathcal{R\cdot G} \\
&\mathcal{=}&\mathcal{R\cdot R},
\end{eqnarray*}%
where $\mathcal{R\cdot C}$ means that $\mathcal{R}\left( X,Y\right) $ acts
as a derivation on $\mathcal{C}$ for any vector fields $X,Y\in \mathfrak{X}%
\left( M\right) $. However,%
\begin{eqnarray*}
\mathcal{C\cdot R} &=&\left( \mathcal{R}-\frac{\tau }{n\left( n-1\right) }%
\mathcal{G}\right) \cdot \mathcal{R} \\
&=&\mathcal{R\cdot \mathcal{R}}-\frac{\tau }{n\left( n-1\right) }\mathcal{G}%
\cdot \mathcal{R}.
\end{eqnarray*}%
We thus have the followings:

\begin{proposition}
A pseudo-Riemannian manifold $M$ admits a semi-symmetric concircular
curvature tensor $\mathcal{C}$ if and only if $M$ is semi-symmetric.
\end{proposition}

\begin{proposition}
A pseudo-Riemannian manifold $M$ is pseudo-symmetric ( i.e. $\mathcal{R\cdot 
\mathcal{R}}=\frac{\tau }{n\left( n-1\right) }\mathcal{G}\cdot \mathcal{R}$)
if and only if $\mathcal{C\cdot R}=0$.
\end{proposition}

On the other hand%
\begin{eqnarray*}
\mathcal{C\cdot C} &=&\left( \mathcal{R}-\frac{\tau }{n\left( n-1\right) }%
\mathcal{G}\right) \cdot \left( \mathcal{R}-\frac{\tau }{n\left( n-1\right) }%
\mathcal{G}\right) \\
&=&\mathcal{R}\cdot \mathcal{\mathcal{R}}-\frac{\tau }{n\left( n-1\right) }%
\mathcal{G}\cdot \mathcal{R}.
\end{eqnarray*}%
We thus have:

\begin{proposition}
A pseudo-Riemannian manifold $M$ is pseudo-symmetric if and only if $%
\mathcal{C\cdot C}=0$.
\end{proposition}

Now, assume that $\mathcal{C}$ vanishes on $M$. Then%
\begin{equation*}
\mathcal{R}=\frac{\tau }{n\left( n-1\right) }\mathcal{G},
\end{equation*}%
i.e., $M$ is of constant curvature $\kappa =\frac{\tau }{n\left( n-1\right) }
$. The converse is also true and we have:

\begin{proposition}
A concircularly flat pseudo-Riemannian manifold $M$ (i.e. $M$ admits a flat
concircular curvature tensor) is of constant curvature.
\end{proposition}

A pseudo-Riemannian manifold $M$ is said to be concircularly symmetric if $%
\nabla \mathcal{C}=0$. It is clear that%
\begin{equation*}
\nabla \mathcal{C}=\nabla \mathcal{R}-\frac{1}{n\left( n-1\right) }\left(
\nabla \tau \right) \mathcal{G}.
\end{equation*}

Assume that $M$ is concircularly symmetric i.e., $\nabla \mathcal{C}=0$. Then%
\begin{equation*}
\nabla \mathcal{R}=\frac{1}{n\left( n-1\right) }\left( \nabla \tau \right) 
\mathcal{G}.
\end{equation*}

The second Bianchi identity implies that $M$ is of constant curvature $%
\kappa $ and consequently $M$ is locally symmetric. Conversely, now suppose
that $M$ is locally symmetric, that is, $\nabla \mathcal{R}=0$, then the
scalar curvature is constant and hence $\nabla \mathcal{C}=0$. This
discussion leads to the following result.

\begin{proposition}
\label{CC-R-1}A pseudo-Riemannian manifold $\left( M,g\right) $ is locally
symmetric if and only if it is concircularly symmetric.
\end{proposition}

In \cite{Boeckx:1993}, it is proved that a semi-symmetric manifold $\left(
M,g\right) $ whose Ricci tensor is a Codazzi tensor is a locally symmetric
manifold. This result and Proposition (\ref{CC-R-1}) lead to the following.

\begin{corollary}
\label{CC-R-2}A semi-symmetric manifold $\left( M,g\right) $ whose Ricci
tensor is a Codazzi tensor is a concircularly symmetric manifold.
\end{corollary}

\section{Concircularly Flat Warped Products}

In this section, we shall first give some basic definitions about warped
product manifolds and then apply them to study the concircularly flat warped
products. Suppose that $\left( M_{1},g_{1},\nabla ^{1}\right) $ and $\left(
M_{2},g_{2},\nabla ^{2}\right) $ are two smooth pseudo-Riemannian manifolds
equipped with Riemannian metrics $g_{i},$ where $\nabla ^{i}$ is the
Levi-Civita connection of the metric $g_{i}$ for $i=1,2.$ Further suppose
that $\pi _{1}\colon M_{1}\times M_{2}\rightarrow M_{1}$ and $\pi _{2}\colon
M_{1}\times M_{2}\rightarrow M_{2}$ are the natural projection maps of the
Cartesian product $M_{1}\times M_{2}$ onto $M_{1}$ and $M_{2},$
respectively. If $f\colon M_{1}\rightarrow \left( 0,\infty \right) $ is a
positive real-valued smooth function, then the warped product manifold $%
M_{1}\times _{f}M_{2}$ is the the product manifold $M_{1}\times M_{2}$
equipped with the metric tensor $g=g_{1}\oplus f^{2}g_{2}$ defined by%
\begin{equation*}
g=\pi _{1}^{\ast }\left( g_{1}\right) \oplus \left( f\circ \pi _{1}\right)
^{2}\pi _{2}^{\ast }\left( g_{2}\right) ,
\end{equation*}%
where $^{\ast }$ denotes the pull-back operator on tensors\cite{Bishop:1969,
Oneill:1983,Shenawy:2016}. The function $f$ is called the warping function
of the warped product manifold $M_{1}\times _{f}M_{2}$. In particular, if $%
f=1$, then $M_{1}\times _{1}M_{2}=M_{1}\times M_{2}$ is the usual Cartesian
product manifold. It is clear that the submanifold $M_{1}\times \{q\}$ is
isometric to $M_{1}$ for every $q\in M_{2}$. Moreover, $\{p\}\times M_{2}$
is homothetic to $M_{2}$. Throughout this article we use the same notation
for a vector field and for its lift to the product manifold\cite%
{El-Sayied:2016,El-Sayied:2017,Shenawy:2015,Shenawy:2016}.

Throughout this section, $\left( M,g,\nabla \right) $ is a (singly) warped
product manifold of $\left( M_{i},g_{i},\nabla ^{i}\right) ,i=1,2$ with
dimensions $n_{i}\neq 1,$ where $n=n_{1}+n_{2}$. $\mathcal{R},\mathcal{R}%
^{i} $ and {\textrm{Ric}}$,${\textrm{Ric}}$^{i}$ denote the curvature tensor
and Ricci curvature tensor on $M,M_{i}$ respectively. Moreover, $\mathrm{grad%
}f,\Delta f$ denote gradient and Laplacian of $f$ on $M_{1}$ and $f^{\sharp
}=f\Delta f+\left( n_{2}-1\right) g_{1}\left( \mathrm{grad}f,\mathrm{grad}%
f\right) $.\ Finally, concircular curvature tensor on $M$ and $M_{i}$ is
denoted by $\mathcal{C}$ and $\mathcal{C}^{i}$ respectively.

We now define generalized Robertson-Walker space-times. Let $(M,g)$ be an $%
n- $dimensional pseudo-Riemannian manifold and $f$ be a positive smooth
function on an open connected subinterval $I$ of $\mathbb{R}$. Then the $%
(n+1)-$dimensional product manifold $I\times M$ furnished with the metric
tensor 
\begin{equation*}
\bar{g}=-\mathrm{d}t^{2}\oplus f^{2}g
\end{equation*}%
is called a generalized Robertson-Walker space-time and is denoted by $\bar{M%
}=I\times _{f}M,$ where $\mathrm{d}t^{2}$ is the Euclidean usual metric
tensor on $I$. These space-times are generalization of the well-known
Robertson-Walker space-times \cite{Sanchez:2000, Sanchez:1998, Sanchez:1999}%
. From now on, we will denote $\frac{\partial }{\partial t}\in \mathfrak{X}%
(I)$ by $\partial _{t}$ to state our results in simpler forms.

Similarly, we define standard static space-times. Let $(M,g)$ be an $n-$%
dimensional pseudo-Riemannian manifold and $f\colon M\rightarrow (0,\infty )$
be a smooth function. Then the $(n+1)-$dimensional product manifold $I\times
M$ furnished with the metric tensor%
\begin{equation*}
\bar{g}=-f^{2}\mathrm{d}t^{2}\oplus g
\end{equation*}%
is called a standard static space-time and is denoted by $\bar{M}%
=I_{f}\times M,$ where $I$ is an open, connected subinterval of $\mathbb{R}$
and $\mathrm{d}t^{2}$ is the Euclidean metric tensor on $I$. Note that
standard static space-times can be considered as a generalization of the
Einstein static universe\cite{Allison:1988, Allison:1998,
Allison:2003,Besse:2008,El-sayied:2017B,El-sayied:2018}.

The following theorem provides a description of the concircular$\ $curvature
tensor on pseudo-Riemannian warped product manifolds.

\begin{proposition}
Let $M=M_{1}\times _{f}M_{2}$ be a singly warped product manifold with the
metric tensor $g=g_{1}\oplus f^{2}g_{2}$. If $X_{i},Y_{i},Z_{i}\in \mathfrak{%
X}(M_{i})$ $i=1,2$, then the concircular curvature tensor $\mathcal{C}$ on $%
M $ is given by%
\begin{eqnarray}
{\mathcal{C}}\left( X_{1},Y_{1}\right) Z_{1} &=&\mathcal{R}^{1}\left(
X_{1},Y_{1}\right) Z_{1}  \notag \\
&&-\frac{\tau }{n\left( n-1\right) }\left[ g_{1}\left( X_{1},Z_{1}\right)
Y_{1}-g_{1}\left( Y_{1},Z_{1}\right) X_{1}\right]  \label{CC-1}
\end{eqnarray}%
\begin{eqnarray}
{\mathcal{C}}\left( X_{2},Y_{1}\right) Z_{1} &=&\left[ \frac{1}{f}%
H^{f}\left( Y_{1},Z_{1}\right) +\frac{\tau }{n\left( n-1\right) }g_{1}\left(
Y_{1},Z_{1}\right) \right] X_{2}  \label{CC-2} \\
{\mathcal{C}}\left( X_{1},Y_{2}\right) Z_{2} &=&fg_{2}\left(
Y_{2},Z_{2}\right) \left[ \nabla _{X_{1}}^{1}\mathrm{grad}f+\frac{\tau f}{%
n\left( n-1\right) }X_{1}\right] ,  \label{CC-3}
\end{eqnarray}%
and%
\begin{eqnarray}
{\mathcal{C}}\left( X_{2},Y_{2}\right) Z_{2} &=&\mathcal{R}^{2}\left(
X_{2},Y_{2}\right) Z_{2}  \notag \\
&&-\left( \left\Vert \mathrm{grad}f\right\Vert _{1}^{2}+\frac{\tau f^{2}}{%
n\left( n-1\right) }\right) \left[ g_{2}\left( X_{2},Z_{2}\right)
Y_{2}-g_{2}\left( Y_{2},Z_{2}\right) X_{2}\right] ,  \label{CC-4}
\end{eqnarray}%
where $H^{f}\left( Y_{1},Z_{1}\right) =g_{1}\left( \nabla _{X_{1}}^{1}%
\mathrm{grad}f,Z_{1}\right) $ is the Hessian of $f$.
\end{proposition}

The following theorem is a direct consequence of the above proposition.

\begin{theorem}
Let $M=M_{1}\times _{f}M_{2}$ be a singly warped product manifold with the
metric tensor $g=g_{1}\oplus f^{2}g_{2}$. $M$ is concircularly flat if and
only if

\begin{enumerate}
\item $M_{1}$ is of constant curvature%
\begin{equation*}
\kappa _{1}=\kappa =\frac{\tau }{n\left( n-1\right) }.
\end{equation*}

\item $\frac{1}{f}H^{f}\left( Y_{1},Z_{1}\right) +\frac{\tau }{n\left(
n-1\right) }g_{1}\left( Y_{1},Z_{1}\right) =0,$ and

\item $M_{2}$ is of constant curvature%
\begin{equation*}
\kappa _{2}=\left\Vert \mathrm{grad}f\right\Vert _{1}^{2}+\frac{\tau f^{2}}{%
n\left( n-1\right) }=\kappa f^{2}+\left\Vert \mathrm{grad}f\right\Vert
_{1}^{2}.
\end{equation*}
\end{enumerate}
\end{theorem}

Now suppose that the concircular curvature tensor $\mathcal{C}$ on $%
M=M_{1}\times _{f}M_{2}$ vanishes, then equation (\ref{CC-2}) implies that 
\begin{equation}
H^{f}\left( Y_{1},Z_{1}\right) =\frac{-\tau f}{n\left( n-1\right) }%
g_{1}\left( Y_{1},Z_{1}\right) ,  \label{CC-5}
\end{equation}%
i.e., $M_{1}$ is of Hessian type. Taking the trace of this equation we get
that%
\begin{equation}
\Delta f=\frac{-n_{1}\tau }{n\left( n-1\right) }f=-n_{1}\kappa _{1}f.
\label{CC-6}
\end{equation}

\begin{corollary}
Let $M=M_{1}\times _{f}M_{2}$ be a concircularly flat singly warped product
manifold with the metric tensor $g=g_{1}\oplus f^{2}g_{2}$. Then $M_{1}$ is
of Hessian type and $\Delta f=-n_{1}\kappa _{1}f$.
\end{corollary}

Now, we note that $\mathcal{C}$ can be simplified if the last position is a
concurrent field. Let $\zeta =\zeta _{1}+\zeta _{2}$ be a vector field on $%
M=M_{1}\times _{f}M_{2}$, $\left\{ e_{i}|1\leq i\leq n_{1}\right\} $ be an
orthonormal basis of $\mathfrak{X}(U_{1})$ and $\left\{ e_{i}|n_{1}+1\leq
i\leq n_{1}+n_{2}\right\} $ be an orthonormal basis of $\mathfrak{X}(U_{2})$
where $U_{i}$ is an open subset of $M_{i}$. Then $\left\{ e_{i}|1\leq i\leq
n_{1}+n_{2}\right\} $ is an orthogonal basis of $\mathfrak{X}\left(
U_{1}\times _{f}U_{2}\right) $. Thus%
\begin{equation*}
\nabla _{e_{i}}\zeta -e_{i}=\nabla _{e_{i}}^{1}\zeta _{1}-e_{i}+e_{i}\left(
\ln f\right) \zeta _{2}
\end{equation*}%
for $1\leq i\leq n_{1}$ and%
\begin{equation*}
\nabla _{e_{i}}\zeta -e_{i}=\zeta _{1}\left( \ln f\right) e_{i}+\nabla
_{e_{i}}^{2}\zeta _{2}-fg_{2}\left( \zeta _{2},e_{i}\right) \mathrm{grad}%
f-e_{i}
\end{equation*}%
for $n_{1}+1\leq i\leq n_{1}+n_{2}$.

\begin{lemma}
Let $M=M_{1}\times _{f}M_{2}$ be a singly warped product manifold with the
metric tensor $g=g_{1}\oplus f^{2}g_{2}$. Then $\zeta =\zeta _{1}+\zeta _{2}$
is a concircular vector field on $M=M_{1}\times _{f}M_{2}$ if and only if $%
\zeta _{1}$ is a concircular vector field on $M_{1}$ and one of the
following conditions holds

\begin{enumerate}
\item $\zeta _{2}$ is a concircular vector field on $M_{2}$, and $f$ is
constant; or

\item $\zeta _{2}=0$ and $\zeta _{1}\left( f\right) =f$.
\end{enumerate}
\end{lemma}

Let $\zeta $ be a concurrent vector field, then%
\begin{equation*}
\mathcal{R}\left( X,Y\right) \zeta =0.
\end{equation*}%
Thus%
\begin{equation*}
\mathcal{C}\left( X,Y\right) \zeta =-\frac{\tau }{n\left( n-1\right) }\left[
g\left( X,\zeta \right) Y-g\left( Y,\zeta \right) X\right] .
\end{equation*}%
Suppose that $M=M_{1}\times _{f}M_{2}$ is a concircularly curvature flat
warped product manifold, then%
\begin{equation*}
\tau \left[ g\left( Y,\zeta \right) X-g\left( X,\zeta \right) Y\right] =0
\end{equation*}%
for any vector fields $X$ and $Y$. Thus $\tau =0$ and consequently $M$ is
flat. This discussion leads to the following result.

\begin{theorem}
Let $M=M_{1}\times _{f}M_{2}$ be a concircularly flat singly warped product
manifold with the metric tensor $g=g_{1}\oplus f^{2}g_{2}$. Then $M$ is flat
if $M_{1}$ admits a concircular vector field $\zeta _{1}$ and one of the
following conditions holds:

\begin{enumerate}
\item $M_{2}$ admits a concircular vector field $\zeta _{2}$ and $f$ is
constant; or

\item $\zeta _{1}\left( f\right) =f$.
\end{enumerate}
\end{theorem}

We will now focus on generalized Robertson-Walker space-times and consider
the concircular curvature on this class of space-times by using our previous
results. Let $\bar{M}=I\times _{f}M$ be a generalized Robertson-Walker
space-time equipped with the metric tensor $\bar{g}=-\mathrm{d}t^{2}\oplus
f^{2}g$. Then the concircular curvature tensor $\mathcal{\bar{C}}$ on $\bar{M%
}$ is given by

\begin{eqnarray*}
\bar{\mathcal{C}}(\partial _{t},\partial _{t})\partial _{t} &=&0,\text{ }{%
\bar{\mathcal{C}}(X,\partial _{t})\partial _{t}=\frac{-1}{f}[\ddot{f}+\dfrac{%
\bar{\tau}f}{n(n+1)}]X}, \\
{\bar{\mathcal{C}}(\partial _{t},X)Y} &{=}&{fg(X,Y)[\ddot{f}+\dfrac{\bar{\tau%
}f}{n(n+1)}]\partial _{t}}, \\
{\bar{\mathcal{C}}(X,Y)Z} &{=}&{\mathcal{R}(X,Y)Z+[\dot{f}^{2}-\dfrac{\bar{%
\tau}f^{2}}{n(n+1)}][g(X,Z)Y-g(Y,Z)X]},
\end{eqnarray*}%
for any vector fields $X,Y,Z\in \mathfrak{X}(M),$ where $\mathcal{R}$ is the
(Riemann) curvature tensor on $M$. By using direct calculation and our
previous results one can conclude the following.

\begin{proposition}
Let $\bar{M}=I\times _{f}M$ be an $\left( n+1\right) -$dimensional
generalized Robertson-Walker space-time equipped with the metric tensor $%
\bar{g}=-\mathrm{d}t^{2}\oplus f^{2}g$, $n\geq 3$. $\bar{M}$ is
concircularly flat if and only if

\begin{enumerate}
\item The scalar curvature of $(\bar{M},\bar{g})$ satisfies ${\ddot{f}+%
\dfrac{\bar{\tau}f}{n(n+1)}=0}$, and

\item $(M,g)$ has constant sectional curvature $\kappa \equiv -\left[ {\dot{f%
}^{2}+f\ddot{f}}\right] $.
\end{enumerate}
\end{proposition}

The above result gives us a full characterization for the warping function $%
f $.

\begin{proposition}
Let $\bar{M}=I\times _{f}M$ be an $\left( n+1\right) -$dimensional
concircularly flat generalized Robertson-Walker space-time equipped with the
metric tensor $\bar{g}=-\mathrm{d}t^{2}\oplus f^{2}g$. Suppose that $\bar{X}%
=h\partial _{t}+X$ is a vector field on $\bar{M},$ where $X$ is a vector
field on $M$ and $h$ is a smooth function on $I.$ Then $(\bar{M},\bar{g})$
is flat if one of the following conditions holds

\begin{enumerate}
\item $M$ admits a concircular vector field and $f$ is constant, or

\item $f(t)=at+b$.
\end{enumerate}
\end{proposition}

Now, we are ready to study concircular curvature tensor $\mathcal{\bar{C}}$
on $\bar{M}=_{f}I\times M$. Let $\bar{M}=I_{f}\times M$ be a standard static
space-time equipped with the metric tensor $\bar{g}=-f^{2}\mathrm{d}%
t^{2}\oplus g$. Then the concircular curvature tensor $\mathcal{\bar{C}}$ on 
$\bar{M}$ is given by%
\begin{eqnarray*}
{{\bar{\mathcal{C}}(X,\partial _{t})}}\partial _{t} &=&{-f}\left[ {{\nabla
_{X}\mathrm{grad}f+{\frac{\bar{\tau}f}{n(n+1)}}X}}\right] , \\
{{\bar{\mathcal{C}}(\partial _{t},X)}}Y &=&\left[ {{{\frac{1}{f}}\mathrm{H}%
^{f}(X,Y)+{\frac{\bar{\tau}}{n(n+1)}}g(X,Y)}}\right] {{\partial _{t}}}, \\
{{\bar{\mathcal{C}}(X,Y}})Z &=&{\mathcal{\bar{R}}(X,Y)Z-{\frac{\bar{\tau}}{%
n(n+1)}}}\left[ {g(X,Z)Y-g(Y,Z)X}\right] ,
\end{eqnarray*}%
for any vector fields $X,Y,Z\in \mathfrak{X}(M),$ where $\mathcal{R}$ is the
(Riemann) curvature tensor on $M$. Now, we can characterize concircularly
flat standard static space-time as:

\begin{proposition}
Let $\bar{M}=I_{f}\times M$ be an $\left( n+1\right) -$dimensional standard
static space-time equipped with the metric tensor $\bar{g}=-f^{2}\mathrm{d}%
t^{2}\oplus g$, $n\geq 3$. $\bar{M}$ is concircularly flat if and only if

\begin{enumerate}
\item {$\nabla $}${{_{X}}}${{$\mathrm{grad}$}}${{f=-\displaystyle{\frac{\bar{%
\tau}}{n(n+1)}}X}}$ for any vector field $X$ on $M$, and

\item $(M,g)$ has constant sectional curvature $\displaystyle{\kappa =\frac{%
\bar{\tau}}{n(n+1)}}$.
\end{enumerate}
\end{proposition}

\section{Concircularly Symmetric Warped Product Manifolds}

A pseudo-Riemannian singly warped product manifold $M$ is said to be
concircular symmetric if%
\begin{equation*}
\left( \nabla _{\zeta }\mathcal{C}\right) \left( X,Y,Z\right) =0
\end{equation*}%
for any vector fields $X,Y,Z$ and $\zeta $. It is clear that (see Section $3$
above)%
\begin{equation*}
\left( \nabla _{\zeta }\mathcal{R}\right) \left( X,Y,Z\right) =0.
\end{equation*}

This condition yields the following consequences%
\begin{equation}
\left( \nabla _{\zeta _{1}}\mathcal{R}\right) \left(
X_{1},Y_{1},Z_{1}\right) =\left( \nabla _{\zeta _{1}}^{1}\mathcal{R}%
^{1}\right) \left( X_{1},Y_{1},Z_{1}\right) =0.  \label{CC-CS-2}
\end{equation}%
Thus $M_{1}$ is locally symmetric. The second case is%
\begin{equation}
\left( \nabla _{\zeta _{1}}\mathcal{R}\right) \left(
X_{2},Y_{1},Z_{1}\right) =0.
\end{equation}%
This yields%
\begin{equation}
-\frac{1}{f^{2}}\zeta _{1}\left( f\right) H^{f}\left( Y_{1},Z_{1}\right)
X_{2}+\frac{1}{f}g_{1}\left( \nabla _{\zeta _{1}}^{1}\nabla _{Y_{1}}^{1}%
\mathrm{grad}f,Z_{1}\right) X_{2}-\frac{1}{f}H^{f}\left( \nabla _{\zeta
_{1}}^{1}Y_{1},Z_{1}\right) X_{2}=0,  \notag
\end{equation}%
i.e. $\mathcal{F}=\frac{1}{f}H^{f}$ is parallel. The next case is%
\begin{equation*}
\left( \nabla _{\zeta _{2}}\mathcal{R}\right) \left(
X_{2},Y_{1},Z_{1}\right) =0
\end{equation*}%
\begin{eqnarray*}
0 &=&\nabla _{\zeta _{2}}\mathcal{R}\left( X_{2},Y_{1}\right) Z_{1}-\mathcal{%
R}\left( \nabla _{\zeta _{2}}X_{2},Y_{1}\right) Z_{1}-Z_{1}\left( \ln
f\right) \mathcal{R}\left( X_{2},Y_{1}\right) \zeta _{2} \\
&=&\mathcal{F}\left( Z_{1},Y_{1}\right) \nabla _{\zeta _{2}}X_{2}-\mathcal{R}%
\left( \nabla _{\zeta _{2}}X_{2},Y_{1}\right) Z_{1}+Z_{1}\left( f\right)
g_{2}\left( X_{2},\zeta _{2}\right) \nabla _{Y_{1}}^{1}\mathrm{grad}f
\end{eqnarray*}%
and so%
\begin{equation*}
\mathcal{R}^{1}\left( \mathrm{grad}f,Y_{1}\right) Z_{1}=\mathcal{F}\left(
Z_{1},Y_{1}\right) \mathrm{grad}f-Z_{1}\left( \ln f\right) \nabla
_{Y_{1}}^{1}\mathrm{grad}f
\end{equation*}%
Now, we have%
\begin{equation*}
\left( \nabla _{\zeta _{1}}\mathcal{R}\right) \left(
X_{2},Y_{2},Z_{2}\right) =\left( \nabla _{\zeta _{2}}\mathcal{R}\right)
\left( X_{1},Y_{2},Z_{2}\right) =0
\end{equation*}%
Thus%
\begin{equation}
X_{1}\left( f\right) \mathcal{R}^{2}\left( \zeta _{2},Y_{2}\right)
Z_{2}=\left( X_{1}\left( f\right) \left\Vert \mathrm{grad}f\right\Vert
^{2}-f^{2}\mathcal{F}\left( X_{1},\mathrm{grad}f\right) \right) G_{2}\left(
\zeta _{2},Y_{2},Z_{2}\right) ,  \label{CC-CS-6}
\end{equation}%
where $G_{2}\left( \zeta _{2},Y_{2},Z_{2}\right) =\left[ g_{2}\left( \zeta
_{2},Z_{2}\right) Y_{2}-g_{2}\left( Y_{2},Z_{2}\right) \zeta _{2}\right] $.
The next case is%
\begin{equation*}
\left( \nabla _{\zeta _{2}}\mathcal{R}\right) \left(
X_{2},Y_{2},Z_{2}\right) =0.
\end{equation*}

This yields%
\begin{equation}
\left( \nabla _{\zeta _{2}}^{2}\mathcal{R}^{2}\right) \left(
X_{2},Y_{2},Z_{2}\right) =0.  \label{CC-CS-8}
\end{equation}

\begin{theorem}
Let $M=M_{1}\times _{f}M_{2}$ be a concircularly symmetric warped product
manifold with the metric tensor $g=g_{1}\oplus f^{2}g_{2}$. Then,

\begin{enumerate}
\item both $M_{1}$ and $M_{2}$ are locally symmetric,

\item $M_{2}$ is of constant curvature given that $f$ is not constant, and

\item $\mathcal{F}=\frac{1}{f}H^{f}$ is parallel.
\end{enumerate}
\end{theorem}

\section{Divergence Free Concircular Curvature Tensor}

It is well known that the Riemann tensor is harmonic if and only if the
Ricci tensor is a Codazzi tensor, i.e. for any vector fields $X,Y,Z\in 
\mathfrak{X}\left( M\right) $, we have%
\begin{equation*}
\left( \nabla _{X}\mathrm{Ric}\right) \left( Y,Z\right) =\left( \nabla _{Y}%
\mathrm{Ric}\right) \left( X,Z\right) .
\end{equation*}%
Moreover, the concircular curvature tensor is divergence free if and only if
the Riemann tensor is harmonic. Let us define%
\begin{equation*}
T\left( X,Y,Z\right) =\left( \nabla _{X}\mathrm{Ric}\right) \left(
Y,Z\right) -\left( D_{Y}\mathrm{Ric}\right) \left( X,Z\right)
\end{equation*}%
for any vector fields $X,Y,Z\in \mathfrak{X}\left( M\right) $. It is clear
that the Ricci tensor is a Codazzi tensor if and only if $T\left(
X,Y,Z\right) $ vanishes. Let $\left( M_{1}\times _{f}M_{2},g\right) $ be a
singly warped product manifold with $T\left( X,Y,Z\right) =0$. Then%
\begin{eqnarray}
T^{1}\left( X_{1},Y_{1},Z_{1}\right) &=&\frac{n_{2}}{f}Y_{1}\left( f\right) 
\mathcal{F}\left( X_{1},Z_{1}\right) -\frac{n_{2}}{f}X_{1}\left( f\right) 
\mathcal{F}\left( Y_{1},Z_{1}\right)  \notag \\
&&-\frac{n_{2}}{f}\mathrm{\mathcal{R}}^{1}\left( X_{1},Y_{1},\mathrm{grad}%
f,Z_{1}\right) .
\end{eqnarray}

The next case is%
\begin{eqnarray*}
0 &=&X_{1}\left( f^{\sharp }\right) g_{2}\left( Y_{2},Z_{2}\right)
-2X_{1}\left( \ln f\right) \mathrm{Ric}\left( Y_{2},Z_{2}\right) \\
&&-Y_{1}\left( f^{\sharp }\right) g_{2}\left( X_{2},Z_{2}\right)
+2Y_{1}\left( \ln f\right) \mathrm{Ric}\left( X_{2},Z_{2}\right)
\end{eqnarray*}

\begin{equation}
X_{1}\left( f\right) \mathrm{Ric}\left( Y_{2},Z_{2}\right) =f\left(
X_{1}\left( f^{\sharp }\right) -f\mathrm{Ric}\left( X_{1},\mathrm{grad}%
f\right) \right) g_{2}\left( Y_{2},Z_{2}\right) .
\end{equation}%
Finally,

\begin{equation}
T^{2}\left( X_{2},Y_{2},Z_{2}\right) =0.
\end{equation}%
The tensor $T$ vanishes in the rest cases. Now, one can write the following
results.

\begin{theorem}
Let $\left( M_{1}\times _{f}M_{2},g\right) $ be a singly warped product
manifold with warping function $f>0$ on $M_{1}$. Assume the concircular
curvature tensor $\mathcal{C}$ is divergence free. Then,

\begin{enumerate}
\item the concircular curvature tensor $\mathcal{C}^{1}$ is divergence free
if%
\begin{equation*}
\mathrm{\mathcal{R}}^{1}\left( X_{1},Y_{1},\mathrm{grad}f,Z_{1}\right)
=Y_{1}\left( f\right) \mathcal{F}\left( X_{1},Z_{1}\right) -X_{1}\left(
f\right) \mathcal{F}\left( Y_{1},Z_{1}\right)
\end{equation*}

\item the concircular curvature tensor $\mathcal{C}^{2}$ is divergence free,
and

\item $f$ is constant or $\left( M_{2},g_{2}\right) $ is Einstein.
\end{enumerate}
\end{theorem}

\begin{theorem}
Let $\left( M_{1}\times _{f}M_{2},g\right) $ be a singly warped product
manifold with warping function $f>0$ on $M_{1}$. The concircular curvature
tensor of the metric tensor $g$ is divergence free if

\begin{enumerate}
\item $f$ is constant and the concircular curvature tensors $\mathcal{C}^{i}$
of the metric tensors $g_{i};i=1,2$ are divergence free, or

\item $H^{f}=0,$ $\mathcal{C}^{1}$ is divergence free and $\left(
M_{2},g_{2}\right) $ is Einstein with factor $g_{1}\left( \mathrm{grad}f,%
\mathrm{grad}f\right) $.
\end{enumerate}
\end{theorem}

The following results are special cases on a generalized Robertson-Walker
space-time and on a standard static space-time.

\begin{corollary}
Let $\bar{M}=I\times _{f}M$ be a generalized Robertson-Walker space-time
with the metric tensor $\bar{g}=-\mathrm{d}t^{2}\oplus f^{2}g$. If the
concircular curvature tensor $\bar{\mathcal{C}}$ of $(\bar{M},\bar{g})$ is
divergence free, then the concircular curvature tensor ${\mathcal{C}}$ of $%
(M,g)$ is divergence free. If, in addition, $f=at+b$, then $(M,g)$ is
Einstein.
\end{corollary}

\begin{corollary}
Let $\bar{M}=I\times _{f}M$ be a generalized Robertson-Walker space-time
with the metric tensor $\bar{g}=-\mathrm{d}t^{2}\oplus f^{2}g.$ Then the
concircular curvature tensor $\bar{\mathcal{C}}$ of $(\bar{M},\bar{g})$ is
divergence free if

\begin{enumerate}
\item $f$ is constant and the concircular curvature tensor ${\mathcal{C}}$
of $(M,g)$ is divergence free, or

\item $f=at+b$ and $(M,g)$ is Einstein with factor $-a{}^{2}$.
\end{enumerate}
\end{corollary}

\begin{corollary}
Let $\bar{M}=_{f}I\times M$ be a standard static space-time with the metric
tensor $\bar{g}=-f^{2}\mathrm{d}t^{2}\oplus g$ and $H^{f}=0$. Then the
concircular curvature tensor $\bar{\mathcal{C}}$ of $(\bar{M},\bar{g})$ is
divergence free if and only if the concircular curvature tensor ${\mathcal{C}%
}$ of $(M,g)$ is divergence free.
\end{corollary}

\section{Acknowledgment}

We would like to thank the referees for their careful reviews and valuable
comments which helped us to improve quality of the paper.


\begin{thebibliography}{99}
\bibitem{Ahsan:2017} Z. Ahsan, \emph{Tensors: Mathematics of Differential
Geometry and Relativity}, PHI Learning Pvt. Ltd., Delhi (Second Printing,
December 2017).

\bibitem{Ahsan:2009} Z. Ahsan and S. A. Siddiqui, \emph{Concircular
Curvature Tensor and Fluid space-times}, Int. J. Theor. Phys., \textbf{48},
3202--3212 (2009).

\bibitem{Allison:1988} D. E. Allison, \emph{Geodesic Completeness in Static
Space-times}, Geometriae Dedicata, \textbf{26}, 85--97 (1988).

\bibitem{Allison:1998} D. E. Allison, \emph{Energy conditions in standard
static space-times}, Gen. Relat. Gravit., \textbf{20}, no.2, 115--122 (1998).

\bibitem{Allison:2003} D. E. Allison and B. \"{U}nal, \emph{Geodesic
Structure of Standard Static Space-times}, J. Geom. Phys., \textbf{46},
no.2, 193--200 (2003).

\bibitem{Besse:2008} A. L. Besse, \emph{Einstein Manifolds}, Classics in
Mathematics, Springer-Verlag, Berlin, (2008).

\bibitem{Bishop:1969} R. L. Bishop and B. O'Neill, \emph{Manifolds of
negative curvature}, Trans. Amer. Math. Soc., \textbf{145}, 1--49 (1969).

\bibitem{Boeckx:1993} E. Boeckx, \emph{Einstein like semi symmetric spaces},
Archivum Mathematicum, Tomus \textbf{29}, 235--240 (1992).

\bibitem{Catino2105} Giovanni Catino, Carlo Mantegazza and Lorenzo Mazzieri 
\emph{A note on Codazzi tensors}, Mathematische Annalen, \textbf{362}, no.
1-2, 629---638 (2015).

\bibitem{Chen2015} Bang-Yen Chen, \emph{Some results on concircular vector
fields and their applications to Ricci solitons}, Bull. Korean Math. Soc. 
\textbf{52}, no. 5, 1535---1547 (2015).

\bibitem{El-Sayied:2016} H. K. El-Sayied, Sameh Shenawy and N. Syied. \emph{%
Conformal vector fields on doubly warped product manifolds and applications}%
, Advances in Mathematical Physics, Volume 2016, Article ID 6508309, 11
pages.

\bibitem{El-Sayied:2017} H. K. El-Sayied, Sameh Shenawy and N. Syied. \emph{%
On symmetries of generalized Robertson-Walker spacetimes and applications},
Journal of Dynamical Systems and Geometric Theories, \textbf{15}, no. 1,
51--69 (2017).

\bibitem{El-sayied:2017B} H. K. El-Sayied, Sameh Shenawy and Noha Syied. 
\emph{Symmetries of }$f-$\emph{associated standard static spacetimes and
applications,} Journal of the Egyptian Mathematical Society 25, no. 4,
414-418 (2017).

\bibitem{El-sayied:2018} H. K. El-Sayied, Sameh Shenawy and Noha Syied. 
\emph{Locally symmetric }$f-$\emph{associated standard static spacetimes,}
Mathematical Methods in the Applied Sciences 41, no. 15, 5733-5736 (2018).

\bibitem{Sanchez:2000} J. L. Flores and M. S\'{a}nchez, \emph{Geodesic
Connectedness and Conjugate Points in GRW Spacetimes}, J. Geom. Phys., 
\textbf{36}, no.3-4, 285--314 (2000).

\bibitem{Gouli:2014} F. Gouli-Andreou and E. Moutafi, \emph{On the
concircular curvature of a }$\left( k,\mu ,\nu \right) $\emph{-manifold},
Pacific Journal of Mathematics, \textbf{269}, no. 1, 113--132 (2014).

\bibitem{Kim1982} In-Bae Kim, \emph{Special concircular vector fields in
Riemannian manifolds}, Hiroshima Math. J. \textbf{12}, no. 1, 77-91 (1982).

\bibitem{Majhi:2013} P. Majhi and U. C. De, \emph{Concircular Curvature
Tensor on }$K-$\emph{Contact Manifolds}, Acta Mathematica Academiae
Paedagogicae Nyregyhaziensis, \textbf{29}, 89--99 (2013).

\bibitem{Majhi:2015} P. Majhi and U. C. De, \emph{Classifications of }$N(k)-$%
\emph{contact metric manifolds satisfying certain curvature conditions, }%
Acta Math. Univ. Comenianae, Vol. LXXXIV, no. 1 , 167--178 (2015).

\bibitem{Oneill:1983} B. O'Neill, \emph{Semi-Riemannian Geometry with
Applications to Relativity}, Academic Press Limited, London, (1983).

\bibitem{Olszak:2012} K. Olszak and Z. Olszak, \emph{On pseudo-Riemannian
manifolds with recurrent concircular curvature tensor}, Acta Mathematica
Hangarica, \textbf{137}, 64-71 (2012).

\bibitem{Pokhariyal:1974} G. P. Pokhariyal, \emph{Curvature tensors in
Riemannian manifolds II}, Proceedings Mathematical Sciences, \textbf{79},
no. 3,105--110 (1974).

\bibitem{Pokhariyal:1982} G. P. Pokhariyal, \emph{Relativistic significance
of curvature tensors}, International Journal of Mathematics and Mathematical
Sciences, \textbf{5}, no.1, 133--139 (1982).

\bibitem{Mishra:1971} G. P. Pokhariyal and R. S. Mishra, \emph{Curvature
tensors in Riemannian manifolds}, Indian Journal of Pure and Applied
Mathematics, \textbf{2}, no. 3, 529--530 (1970).

\bibitem{Pokhariyal:1970} G. P. Pokhariyal and R. S. Mishra, \emph{Curvature
tensors and their relativistics significance}, Yokohama Math. J., \textbf{18}%
, 105--108 (1970).

\bibitem{Sanchez:1999} M. S\'{a}nchez, \emph{On the Geometry of Generalized
Robertson-Walker spacetimes: Curvature and Killing fields}, J. Geom. Phys., 
\textbf{31}, no.1, 1--15 (1999).

\bibitem{Sanchez:1998} M. S\'{a}nchez, \emph{\ On the Geometry of
Generalized Robertson-Walker spacetimes: geodesics}, Gen. Relat. Gravit., 
\textbf{30}, no.6, 915--932 (1998).

\bibitem{Singh:1999} H. Singh and Q. Khan, \emph{On symmetric Riemannian
manifolds}, Novi Sad Journal of Mathematics, \textbf{29}, no. 3, 301--308
(1999).

\bibitem{Shenawy:2015} Sameh Shenawy and B\"{u}lent \"{U}nal. $2-$\emph{%
Killing vector fields on warped product manifolds}, International Journal of
Mathematics, \textbf{26}(2015), 17 pages.

\bibitem{Shenawy:2016} Sameh Shenawy and B\"{u}lent \"{U}nal. \emph{The }$%
W_{2}$\emph{-curvature tensor on warped product manifolds and applications},
International Journal of Geometric Methods in Modern Physics, \textbf{13},
no. 7 (2016) 1650099 (14 pages).

\bibitem{Yano:1940} K. Yano, \emph{Concircular geometry I. Concircular
transformations}, Proceedings of the Imperial Academy 16, no. 6 195-200
(1940).

\bibitem{Yano:1984} K. Yano and M. Kon, \emph{Structures on Manifolds},
World Scientific Publishing, Singapore, (1984).

\bibitem{Youssef:2013} N. L. Youssef and A. Soleiman, \emph{On concircularly
recurrent Finsler manifolds}, Balk. J. Geom. Appl., \textbf{18}, no. 1,
101--113 \textbf{(}2013).

\bibitem{Zlatanovic:2014} M. Zlatanovica, I. Hinterleitnerb and M.
Najdanovi, \emph{On Equitorsion Concircular Tensors of Generalized
Riemannian Spaces}, Filomat \textbf{28}, no. 3, 463--471 (2014).
\end{thebibliography}
\end{document}